\documentclass[reqno]{amsart}
\usepackage{amssymb}
\usepackage{amsthm}
\usepackage{amsmath}
\usepackage{a4wide}
\usepackage{latexsym}
\usepackage[v2,tips]{xy}
\usepackage{mathrsfs}
\usepackage[T1]{fontenc}
\usepackage[latin1]{inputenc}
\usepackage{enumitem}
\usepackage{graphicx}
\usepackage{hyperref}
\usepackage[normalem]{ulem}
\usepackage{dutchcal}
\usepackage{cleveref}
\usepackage{tabularx}
\usepackage{color}

\newcommand{\N}{\mathbb{N}}

\newcommand{\R}{\mathbb{R}}
\newcommand{\C}{\mathbb{C}}
\newcommand{\K}{\mathbb{K}}

\newtheorem{theorem}{Theorem}[section]
\newtheorem{lemma}[theorem]{Lemma}
\newtheorem{proposition}[theorem]{Proposition}
\newtheorem{corollary}[theorem]{Corollary}
\theoremstyle{definition}
\newtheorem{definition}[theorem]{Definition}

\newtheorem{remark}[theorem]{Remark}

\newtheorem*{acknowledgements}{Acknowledgements}

\numberwithin{equation}{section}

\renewcommand{\ge}{\geqslant}

\renewcommand{\le}{\leqslant}
\renewcommand{\leq}{\leqslant}

\newcommand{\SC}{\operatorname{SC}}

\title[Compact compositions of strictly singular operators]{Compactness of compositions of\\ strictly singular operators on direct sums of\\ Baernstein, Schreier and $\ell_p$-spaces}
\author[N.J.~Laustsen]{Niels Jakob Laustsen}
\address{(N.J.~Laustsen) School of Mathematical Sciences, Fylde
  College, Lancaster University, Lancaster LA1 4YF, United Kingdom}
\email{n.laustsen@lancaster.ac.uk}

\author[H.~Wirzenius]{Henrik Wirzenius}
\address{(H.~Wirzenius) Institute of Mathematics, Czech Academy of Sciences, \v{Z}itn\'a 25, 115 67 Praha 1, Czech Republic}
\email{wirzenius@math.cas.cz, henrik.wirzenius@gmail.com}


\begin{document}

\begin{abstract} Let~$X$ be the direct sum of finitely many Banach spaces chosen from the following three families: (i)~the Baernstein spaces~$B_p$ for $1<p<\infty$; (ii)~the $p$-convexified Schreier spaces~$S_p$ for $1\le p<\infty$; (iii)~the sequence spaces~$\ell_p$ for $1\le p<\infty$ (and~$c_0$). We show that the quotient algebra of strictly singular by compact operators on~$X$ is nilpotent; that is, there is a natural number~$k$, dependent only on the collections of direct summands from each of the three families, such that: 
\begin{itemize}
\item every composition of~$k+1$ strictly singular operators on~$X$ is compact; 
\item there are~$k$ strictly singular operators on~$X$ whose composition is not compact.
\end{itemize}
\smallskip

\noindent Accepted to appear in \emph{Proceedings of the American Mathematical Society.}
\end{abstract}

\maketitle

\section{Introduction}

\noindent 
Ever since Kato~\cite{K} introduced the class of strictly singular operators nearly 70 years ago as a larger ideal that can play the same role as the ideal of compact operators in the perturbation theory of Fredholm operators, it has been an important task to describe the precise relationship between the strictly singular and compact operators on a given Banach space. 

The simplest scenario is that these two ideals are equal; this happens for instance for the classical sequence spaces~$\ell_p$, $1\le p<\infty$, and~$c_0$, the quasi-reflexive James spaces~$J_p$, $1<p<\infty$, and the Tsirelson space~$T$, but is quite rare. 

A more common phenomenon is that every composition of two strictly singular operators is compact; Banach spaces with this property include the Lebesgue spaces~$L_p[0,1]$ for $1\le p<\infty$, the continuous functions~$C(K)$ on a compact Hausdorff space~$K$, the direct sums~$\ell_p\oplus\ell_q$ and~$\ell_p\oplus c_0$ for $1\le p<q<\infty$, and---very importantly for us---the $p^{\text{th}}$ Baernstein space~$B_p$ for $1<p<\infty$ and the $p$-convexified Schreier space~$S_p$ for $1\le p<\infty$ (see \cite[Theorem~1.1]{LS2}; we refer to~\Cref{S:2} for the precise definitions of the spaces~$B_p$ and~$S_p$).

There is nothing special about the number~$2$ here, of course; for any $k\in\N$, there are examples of Banach spaces~$X$ with the property that every composition of~$k+1$ strictly singular operators on~$X$ is compact, but there are~$k$ strictly singular operators whose composition is not compact. In algebraic parlance, this means that the quotient algebra of strictly singular by compact operators on~$X$ is nilpotent of index~$k+1$.  A recent example where this occurs is the Schreier space $X[\mathcal{S}_k]$ induced by the $k^{\text{th}}$ Schreier family~$\mathcal{S}_k$ (see \cite[Theorem~7.6(3)--(4)]{CPB}; note that this theorem also applies to the Schreier space~$X[\mathcal{S}_\xi]$ induced by the Schreier family~$\mathcal{S}_\xi$ for a countably infinite ordinal~$\xi$, but the value of the index of nilpotency is harder to state explicitly in those cases). Many other examples exist, including Tarbard's variant~$\mathfrak{X}_{k+1}$ of the Argyros--Haydon space~\cite{T} and  the Tsirelson-like space~$\mathfrak{X}_{0,1}^k$ constructed by Argyros, Beanland and Motakis~\cite{ABP}.

Our main result builds on several of the above results, as we consider nilpotency of the quotient algebra of strictly singular by compact operators on Banach spaces that are finite direct sums of Baernstein spaces, $\ell_p$-spaces and $p$-convexified Schreier spaces. The precise statement is as follows. 

\begin{theorem}\label{Thm:main}
    Let $L\subset(1,\infty)$ and $M,N\subset[1,\infty)$ be finite sets, not all empty, and define
    \begin{equation}\label{Thm:main:eqk} k = \begin{cases} \lvert L\rvert+\lvert L\cup M\rvert -1\ &\text{if}\ N=\emptyset\\ \lvert L\rvert + \lvert L\cup M\rvert + \lvert N\rvert\ &\text{otherwise.} \end{cases} \end{equation}    
    Then the Banach space
    \begin{equation}\label{Thm:main:eqX} X = \Bigl(\bigoplus_{p\in L} B_p\Bigr)\oplus\Bigl(\bigoplus_{q\in M}\ell_q\Bigr)\oplus\Bigl(\bigoplus_{r\in N} S_r\Bigr) \end{equation}
    satisfies:
    \begin{enumerate}[label={\normalfont{(\roman*)}}]  
    \item\label{thm:main:item1} every composition of~$k+1$ strictly singular operators on~$X$ is compact; 
    \item\label{thm:main:item2} there are~$k$ strictly singular operators on~$X$ whose composition is not compact.
    \end{enumerate}
    Hence, the quotient algebra $\mathscr{S}(X)/\mathscr{K}(X)$ of strictly singular by compact operators on~$X$ is nilpotent of index~$k+1$. 
\end{theorem}

\begin{remark}\label{R:LMdisjoint}
\begin{enumerate}[label={\normalfont{(\roman*)}}]  
    \item\label{R:LMdisjoint1} The statement of \Cref{Thm:main} does not specify which norm we put on the (finite) direct sum~\eqref{Thm:main:eqX} defining the Banach space~$X$ because the conclusions depend only on the isomorphism class of~$X$. For definiteness, we may choose the norm $\lVert (x_j)_{j=1}^n\rVert_\infty = \max_{1\le j\le n}\lVert x_j\rVert$, where $n = \lvert L\rvert + \lvert M\rvert + \lvert N\rvert$, but emphasize that the theorem is true for any equivalent norm such as $\lVert (x_j)_{j=1}^n\rVert_p = \bigl(\sum_{j=1}^n\lVert x_j\rVert^p\bigr)^{1/p}$ for $1\le p<\infty$.
\item\label{R:LMdisjoint1.5} Important special cases of \Cref{Thm:main} occur when two of the three index sets~$L$, $M$ and~$N$ are empty. This result is well-known for $X = \bigoplus_{j=1}^m\ell_{p_j}$, where \mbox{$1\le p_1<p_2<\cdots<p_m<\infty$} (see~\cite[Theorem~4.7]{TW} for details), but the other two cases appear to be new for $m\ge 2$, so we state them explicitly for future reference: 
\begin{itemize} 
\item Let $X= \bigoplus_{j=1}^m B_{p_j}$ for some $m\in\N$ and $1< p_1<p_2<\cdots<p_m<\infty$. Then the quotient algebra $\mathscr{S}(X)/\mathscr{K}(X)$ is nilpotent of index~$2m$. 
\item Let $X=\bigoplus_{j=1}^m S_{p_j}$ for some $m\in\N$ and $1\le p_1<p_2<\cdots<p_m<\infty$. Then the quotient algebra $\mathscr{S}(X)/\mathscr{K}(X)$ is nilpotent of index~$m+1$. 
\end{itemize}    
(As already mentioned, \cite[Theorem~1.1]{LS2} contains these results for $m=1$; see also \Cref{L:compactnessofSScompositions} below.)
    \item\label{R:LMdisjoint2} The reason that the cardinality of the set~$L\cup M$  appears in the formula~\eqref{Thm:main:eqk} for~$k$ is that on the one hand,  $B_p$ contains a complemented copy of~$\ell_p$ for $1<p<\infty$, so~$X$ contains a complemented copy of~$\ell_p$ for every  $p\in L\cup M$, and on the other, $\ell_p\oplus\ell_p\cong\ell_p$, so those~$p$ that belong to~$L\cap M$ contribute only one copy of~$\ell_p$ to the collection of complemented subspaces of~$X$. 
    
    Building on these arguments, we see that since $B_p\cong W\oplus\ell_p$ for some Banach space~$W$,  we have $B_p\oplus\ell_p\cong W\oplus\ell_p\oplus\ell_p\cong W\oplus\ell_p\cong B_p$, and therefore $X\cong X\oplus\ell_p$ for $p\in L$. Consequently, we can replace the index set~$M$ in the definition~\eqref{Thm:main:eqX} of~$X$ with~$M\setminus L$ or~$L\cup M$, or any set between them, without affecting the isomorphism class of~$X$.
\item\label{R:LMdisjoint:4} As indicated in the abstract, there is a variant of \Cref{Thm:main} that includes~$c_0$ in the direct sum. To state it concisely, set $Y=X\oplus c_0$, where~$X$ is the Banach space defined by~\eqref{Thm:main:eqX}. Then:
\begin{itemize}
    \item $Y$ is isomorphic to~$X$ if $N\ne\emptyset$, so the conclusions of \Cref{Thm:main} apply verbatim to~$Y$ in this case.
    \item Otherwise $Y = \bigl(\bigoplus_{p\in L} B_p\bigr)\oplus\bigl(\bigoplus_{q\in M}\ell_q\bigr)\oplus c_0$ for some finite sets $L\subset(1,\infty)$ and $M\subset[1,\infty)$, and we have: Every composition of $\lvert L\rvert + \lvert L\cup M\rvert+1$ strictly singular operators on~$Y$ is compact, but there are $\lvert L\rvert + \lvert L\cup M\rvert$ strictly singular operators on~$Y$ whose composition is not compact.
\end{itemize}
We refer to \Cref{R:includingc0} for a detailed justification of these two claims.
\end{enumerate}
\end{remark}

\section{Preliminaries}\label{S:2}

\noindent All vector spaces (in particular Banach spaces) are over the same scalar field~$\K$, which is either~$\R$ or~$\C$. We use function notation for sequences, thus writing $x(n)$ for the $n^{\text{th}}$ coordinate of a sequence~$x\in\K^\N$.  As usual, $c_{00}$ denotes the subspace of finitely supported elements of~$\K^\N$, and 
$(e_{n})_{n\in\N}$ is the \emph{unit vector basis} given by $e_n(m) = 1$ if $m=n$ and $e_n(m) = 0$ otherwise. 

By an \emph{operator,} we mean a bounded linear map between Banach spaces, and write~$\mathscr{B}(X,Y)$ for the set of operators from a Banach space~$X$ to a Banach space~$Y$, abbreviated~$\mathscr{B}(X)$ when $X=Y$; $I_X$ denotes the identity operator on~$X$. An operator~$T\in\mathscr{B}(X,Y)$ is \emph{strictly singular} if no restriction of~$T$ to an in\-finite-di\-men\-sional subspace of~$X$ is an isomorphic embedding. We write~$\mathscr{S}(X,Y)$ and~$\mathscr{K}(X,Y)$ for the sets of strictly singular and compact operators from~$X$ to~$Y$, respectively, abbreviated~$\mathscr{S}(X)$ and $\mathscr{K}(X)$  when $X=Y$.  They are closed operator ideals in the sense of Pietsch, and $\mathscr{K}(X,Y)\subseteq\mathscr{S}(X,Y)$ for any  Banach spaces~$X$ and~$Y$. 

Before we can introduce our main objects of interest---the Schreier and Baernstein spaces---we require the following notion, originating in~\cite{schreier}:  A~\emph{Schreier set} is a finite subset~$F$ of the natural numbers $\N=\{1,2,3,\ldots\}$ such that either $F=\emptyset$ or $\lvert F\rvert\leq \min F$, where $\lvert F\rvert$ denotes the cardinality of~$F$. As usual, we write~$\mathcal{S}_{1}$ for the family of Schreier sets.     
Following~\cite[Section~3]{BL}, for $1\le p<\infty$ and $x\in\K^\N$, we define
\[ \lVert x\rVert_{S_{p}} = \sup\{ \mu_{p}(x,F) : F \in \mathcal{S}_{1}\}\in[0,\infty],\quad\text{where}\quad   \mu_{p}(x,F) =\begin{cases} 0\ &\text{if}\ F=\emptyset,\\ \displaystyle{\Bigl(\sum_{n \in F}\lvert x(n)\rvert^{p}\Bigr)^{1/p}}\ &\text{otherwise.} \end{cases} \] 
Then $\mu_p(\,\cdot\,,F)$ is a seminorm on~$\K^\N$, and $Z_p =\{ x\in\K^\N : \lVert x \rVert_{S_{p}} < \infty\}$ is a subspace of~$\K^\N$ on which~$\lVert\,\cdot\,\rVert_{S_p}$ defines a complete norm. However, the Banach space~$Z_p$ fails to be separable (see \cite[Corollary~5.6]{BL}), so we define the \emph{$p$-convexified Schreier space}, denoted~$S_p$, as the closure of~$c_{00}$ in~$Z_p$. The unit vector basis $(e_n)_{n\in\N}$ is a $1$\nobreakdash-un\-con\-di\-tional, shrinking, normalized basis for~$S_p$ by \cite[Propo\-si\-tions~3.5 and~3.10 and Corollary~3.12]{BL}. 

The analogous definition of the Baernstein spaces involves the notion of a \emph{Schreier chain}, which is a non-empty, finite collection~$\mathcal{C}$ of non-empty, consecutive Schreier sets; that is, 
$\mathcal{C} = \{F_{1},\ldots,F_{m}\}$ for some $m\in\N$ and $F_1,\ldots,F_m\in\mathcal{S}_1\setminus\{\emptyset\}$ with $\max F_{j} < \min F_{j+1}$ for $1\le j <m$. 
Writing~$\SC$ for the collection of all Schreier chains, for $1<p<\infty$ and $x\in\K^\N$, we can define
\[ \lVert x\rVert_{B_{p}} = \sup\{\beta_{p}(x,\mathcal{C}) : \mathcal{C}\in\SC\}\in[0,\infty], \quad\text{where}\quad
\beta_{p}(x,\mathcal{C}) = \biggl(\sum_{F\in \mathcal{C}}\Bigl(\sum_{n\in F}\lvert x(n)\rvert\Bigr)^{p}\biggr)^{1/p}. \]
As before, $\beta_p(\,\cdot\,,\mathcal{C})$ is a seminorm on~$\K^\N$, and $B_p =\{ x\in\K^\N : \lVert x \rVert_{B_{p}} < \infty\}$ is a subspace of~$\K^\N$ on which~$\lVert\,\cdot\,\rVert_{B_p}$ defines a complete norm. In contrast to the Schreier spaces, $c_{00}$ is dense in~$B_p$, which is the $p^{\text{th}}$ \emph{Baernstein space}. It is reflexive, and the unit vector basis is a $1$\nobreakdash-un\-con\-di\-tional, normalized basis for it. Baernstein~\cite{B} originally defined~$B_2$, while Seifert~\cite{S} observed that Baernstein's definition works for general~$p>1$. 

We conclude this preliminary section with a standard piece of terminology from algebra that we already used in the statement of \Cref{Thm:main}: A ring~$\mathscr{A}$ is \emph{nilpotent} if, for some $k\in\N$, we have $a_1a_2\cdots a_k=0$ whenever $a_1,\ldots,a_k\in\mathscr{A}$; the smallest value of~$k$ for which this identity is satisfied is called the \emph{index of nilpotency} of~$\mathscr{A}$.
\section{The proof of Theorem~\ref{Thm:main}}
\noindent 
We present the proofs of the two parts of the theorem separately, beginning with the first. It relies on the matrix representation of operators on a finite direct sum of Banach spaces, defined as follows. 

\begin{definition} Let $X = \bigoplus_{j=1}^n X_j$ for some Banach spaces $X_1,\ldots,X_n$. The \emph{matrix} associated with an operator $T\in\mathscr{B}(X)$ is the operator-valued $(n\times n)$-matrix $(T_{j,k})_{j,k=1}^n$ given by 
\[ T_{j,k}= Q_jTJ_k\in\mathscr{B}(X_k,X_j)\qquad (j,k\in \{1,\ldots,n\}), \] where $Q_j\colon X\to X_j$ and $J_k\colon X_k\to X$ denote the $j^{\text{th}}$ coordinate projection and $k^{\text{th}}$ coordinate embedding, respectively. 
\end{definition}
It is easy to see that composition of operators corresponds to matrix multiplication:
\begin{equation}\label{Eq:matrixmult} 
(TU)_{j,m} = \sum_{k=1}^n T_{j,k}U_{k,m}\qquad (T,U\in\mathscr{B}(X),\,j,m\in \{1,\ldots,n\}). \end{equation}
Furthermore, the identity $T=\sum_{j,k=1}^n J_j T_{j,k} Q_k$ implies that, for any operator ideal~$\mathscr{I}$, we have
\begin{equation}\label{Eq:opidealmatrix} T\in\mathscr{I}(X)\quad \iff\quad T_{j,k}\in\mathscr{I}(X_k,X_j)\ \text{for every}\ j,k\in\{1,\ldots,n\}. \end{equation}

We shall also require the following result,  most of which is known.
\begin{lemma}\label{L:compactnessofSScompositions}  \begin{enumerate}[label={\normalfont{(\roman*)}}]  
    \item\label{L:compactnessofSScompositions1} Every strictly singular operator on~$\ell_p$ is compact for $1\le p<\infty$.
    \item\label{L:compactnessofSScompositions2} Every composition of two strictly singular operators on~$B_p$ is compact for $1< p<\infty$.
    \item\label{L:compactnessofSScompositions3} The composite operator $TU$ is compact whenever $U\colon S_p\to S_p$ and $T\colon S_p\to Y$ are strictly singular, where $1\le p<\infty$ and~$Y$ can be any Banach space. 
    \end{enumerate}
\end{lemma}
\begin{proof}
    Part~\ref{L:compactnessofSScompositions1} is well known; we refer to the sentence after the proof of \cite[Proposition~2.c.3]{LT1} for details. 
    Part~\ref{L:compactnessofSScompositions2} is due to Laustsen and Smith \cite[Theorem~1.1]{LS2}; so is part~\ref{L:compactnessofSScompositions3}, but only for $Y=S_p$. However, invoking a theorem of Rosenthal, we can easily extend their proof to the general case. Indeed, suppose that $U\in\mathscr{S}(S_p)$ and $T\in\mathscr{B}(S_p,Y)$ are operators whose composition~$TU$ is not compact. Our aim is to prove that~$T$ is not strictly singular. Using the elementary observation stated in \cite[Lemma~3.4]{LS2}, we can find a normalized block basic sequence~$(u_n)_{n\in\N}$ of the unit vector basis for~$S_p$ such that \mbox{$\inf_{n\in\N}\lVert TUu_n\rVert_Y>0$.} Then also $\inf_{n\in\N}\lVert Uu_n\rVert_{S_p}>0$, and $(u_n)_{n\in\N}$ is weakly null because the unit vector basis for~$S_p$ is shrinking, so \cite[Lemma~3.3(ii)]{LS2} implies that  $(u_n)_{n\in\N}$ admits a subsequence $(u_{n_j})_{j\in\N}$ such that $(Uu_{n_j})_{j\in\N}$ is equivalent to the unit vector basis $(d_j)_{j\in\N}$ for~$c_0$; let $V\in\mathscr{B}(c_0,S_p)$ be the operator given by $Vd_j= Uu_{n_j}$ for $j\in\N$. Since \[ \inf_{j\in\N}\lVert TVd_j\rVert_Y = \inf_{j\in\N}\lVert TUu_{n_j}\rVert_Y>0, \]  a famous result of Rosenthal, originally stated as the first remark following \cite[Theorem~3.4]{ROS}, implies that~$\N$ contains an infinite subset~$N$ for which the restriction of~$TV$ to the closed span of \mbox{$\{ d_j : j\in N\}$} is an isomorphic embedding. Hence~$T$ is not strictly singular.  
\end{proof}

\begin{proof}[Proof of Theorem~{\normalfont{\ref{Thm:main}}}\ref{thm:main:item1}] In view of \Cref{R:LMdisjoint}\ref{R:LMdisjoint2}, we may suppose that the sets~$L$ and~$M$ are disjoint.  
Our aim is to prove that the composition of~$k+1$ strictly singular operators $R^{(1)},\ldots,R^{(k+1)}$ on~$X$ is compact, so by~\eqref{Eq:opidealmatrix}, we must show that the $(j,m)^{\text{th}}$ entry of the matrix of the composite operator $R^{(k+1)} R^{(k)}\cdots R^{(1)}$ is compact for every $j,m\in\{1,\ldots,n\}$, where $n = \lvert L\rvert+\lvert M\rvert+\lvert N\rvert$. Applying the identity~\eqref{Eq:matrixmult} repeatedly, we obtain
\[ (R^{(k+1)} R^{(k)}\cdots R^{(1)})_{j,m} = \sum_{i_1,\ldots,i_k=1}^n R^{(k+1)}_{j,i_k} R^{(k)}_{i_k,i_{k-1}}\cdots R^{(1)}_{i_1,m}. \]

We shall now complete the proof by showing that each of the $n^k$ terms on the right-hand side of this equation is compact.  Simplifying the notation, we see that this amounts to verifying 
that the composite operator $T:=T_{k+1}T_k\cdots T_1$ is compact whenever $T_j\colon X_j\to X_{j+1}$ is a strictly singular operator for every $1\le j\le k+1$ and the Banach spaces $X_1,\ldots,X_{k+2}$ belong to the family
$\{ B_p : p\in L\}\cup\{\ell_q : q\in M\}\cup \{S_r : r\in N\}$. 

For integers $1\le i\le j\le k+2$, set 
\[ T_{(i\to j)} = \begin{cases} I_{X_j}\ &\text{if}\ i=j,\\ T_{j-1}\cdots T_{i+1}T_i\in\mathscr{S}(X_i,X_j)\ &\text{otherwise.} \end{cases} \] 
This  notation will allow us to justify the following three observations concisely: 
\begin{enumerate}[label={\normalfont{(\roman*)}}] 
    \item\label{case2} Suppose that $X_h=X_i=X_j=B_p$ for some $p\in L$  and integers $1\le h<i<j\le k+2$. Then $T_{(h\to i)},T_{(i\to j)}\in\mathscr{S}(B_p)$, so their composition is compact by \Cref{L:compactnessofSScompositions}\ref{L:compactnessofSScompositions2}, and therefore $T=T_{(j\to k+2)}T_{(i\to j)} T_{(h\to i)}T_{(1\to h)}$ is compact. 
    \item\label{case1} Suppose that $X_i=X_j=\ell_q$ for some $q\in M$ and integers $1\le i<j\le k+2$. Then $T_{(i\to j)}\in\mathscr{S}(\ell_q)$ is compact by \Cref{L:compactnessofSScompositions}\ref{L:compactnessofSScompositions1}, so $T=T_{(j\to k+2)} T_{(i\to j)}T_{(1\to i)}$ is compact.    
    \item\label{case3} Suppose that $X_i=X_j=S_r$ for some $r\in N$  and integers $1\le i<j\le k+1$. Then $T_{(i\to j)}\in\mathscr{S}(S_r)$ and $T_{(j\to k+2)}\in\mathscr{S}(S_r,X_{k+2})$ because $j<k+2$, so \Cref{L:compactnessofSScompositions}\ref{L:compactnessofSScompositions3} implies that their composition is compact. Hence $T=T_{(j\to k+2)}T_{(i\to j)}T_{(1\to i)}$ is compact. 
\end{enumerate}   

If $N=\emptyset$, then $k = 2\lvert L\rvert+\lvert M\rvert-1$ because $L\cap M=\emptyset$. Therefore, choosing $k+2$ spaces $X_1,\ldots,X_{k+2}$ from the family $\{ B_p : p\in L\}\cup\{\ell_q : q\in M\}$, we are either in case~\ref{case2} or~\ref{case1}, so~$T$ is compact. 

Otherwise $N\ne\emptyset$, and we have $k = 2\lvert L\rvert + \lvert M\rvert + \lvert N\rvert$. If we are not in cases~\ref{case2} or~\ref{case1}, then at least $\lvert N\rvert + 2$ of the spaces $X_1,\ldots,X_{k+2}$ come from the family~$\{S_r : r\in N\}$, so we must be in case~\ref{case3}; hence~$T$ is compact. 
\end{proof}

In order to streamline the presentation of the proof of the second part of \Cref{Thm:main}, we require a few preparations. 

\begin{lemma}\label{L:complemented}
    Let $X$ be a Banach space, and~$\mathscr{I}$ and~$\mathscr{J}$ operator ideals. Suppose that, for some $k\in\N$, $X$~contains $k+1$ complemented subspaces $X_1,\ldots,X_{k+1}$ for which there are~$k$ operators $R_1\in\mathscr{J}(X_1,X_2),\ldots,R_k\in\mathscr{J}(X_k,X_{k+1})$ whose composition $R_kR_{k-1}\cdots R_1$ does not belong to $\mathscr{I}(X_1,X_{k+1})$. Then there are operators $T_1,\ldots,T_k\in\mathscr{J}(X)$ whose composition $T_kT_{k-1}\cdots T_1$ does not belong to $\mathscr{I}(X)$. 
\end{lemma}

\begin{proof}
    For each $1\le j\le k+1$, we can take operators $U_j\in\mathscr{B}(X,X_j)$ and $V_j\in\mathscr{B}(X_j,X)$ such that $U_jV_j=I_{X_j}$ because~$X_j$ is complemented in~$X$. Set $T_j = V_{j+1}R_jU_j\in\mathscr{J}(X)$ for $1\le j\le k$. Then we have
    \begin{align*} U_{k+1}(T_kT_{k-1}\cdots T_1)V_1 &= (U_{k+1} V_{k+1})R_k (U_kV_k) R_{k-1}(U_{k-1}V_{k-1})\cdots (U_2V_2)R_1(U_1V_1)\\ &= R_k R_{k-1}\cdots R_1\notin \mathscr{I}(X_1,X_{k+1}), \end{align*}
    so $T_kT_{k-1}\cdots T_1\notin \mathscr{I}(X)$ because~$\mathscr{I}$ is an operator ideal. 
\end{proof}

``Formal inclusion maps'' are at the heart of our proof of \Cref{Thm:main}\ref{thm:main:item2}, so our next step is to make this notion precise; we use a definition specifically tailored to the context at hand, where inclusion between Banach spaces is unambiguous because all the Banach spaces we consider consist of scalar sequences, equipped with the coordinatewise vector space operations inherited from~$\K^\N$. 

\begin{definition}\label{Def:FIM}
    Let $Y$ and $Z$ be vector subspaces of~$\K^\N$, equipped with norms~$\lVert\,\cdot\,\rVert_Y$ and~$\lVert\,\cdot\,\rVert_Z$, respectively. We say that the pair $(Y,Z)$ \emph{admits a formal inclusion map} if $Y\subseteq Z$ and there is a constant~$C_{Y,Z}>0$ such that 
   \begin{equation}\label{Def:FIM:eq} \lVert y\rVert_Z\le C_{Y,Z}\lVert y\rVert_Y\qquad (y\in Y).
   \end{equation}   
\end{definition}

Obviously, the significance of this definition is that when $(Y,Z)$ admits a formal inclusion map, the restriction of the identity operator on~$\K^\N$ defines a bounded linear map from~$Y$ to~$Z$ with norm at most~$C_{Y,Z}$; we denote this map $R_{Y,Z}\colon Y\to Z$ and call it the \emph{formal inclusion map} from~$Y$ to~$Z$, and remark that~$R_{Y,Z}$ is simply the identity operator on~$Y$ if $Y=Z$.

\begin{lemma}\label{L:c00}
  Let~$Y$ and~$Z$ be vector subspaces of~$\K^\N$, equipped with complete norms~$\lVert\,\cdot\,\rVert_Y$ and~$\lVert\,\cdot\,\rVert_Z$, respectively, and suppose that the unit vector basis $(e_n)_{n\in\N}$ is a Schauder basis for both~$Y$ and~$Z$.  Then $(Y,Z)$ admits a formal inclusion map if (and only if) there is a con\-stant~$C_{Y,Z}>0$ such that 
  \begin{equation}\label{L:c00:eq}
      \lVert x\rVert_Z\le C_{Y,Z}\lVert x\rVert_Y\qquad (x\in c_{00}). 
  \end{equation}
\end{lemma}  

\begin{proof}  To prove the non-trivial implication, suppose that~\eqref{L:c00:eq} is satisfied, and take $y\in Y$. The fact that~$(e_n)_{n\in\N}$ is a basis for~$Y$ means that $\lVert y-P_ny\rVert_Y\to 0$ as $n\to\infty$, where $P_n\colon Y\to c_{00}$ denotes the $n^{\text{th}}$ basis projection given by $P_ny=\sum_{j=1}^n y(j)e_j$ for $n\in\N$. In particular, $(P_ny)_{n\in\N}$  is a Cauchy sequence in~$c_{00}$ with respect to the norm~$\lVert\,\cdot\,\rVert_Y$ and therefore also with respect to~$\lVert\,\cdot\,\rVert_Z$ by~\eqref{L:c00:eq}, so $(P_ny)_{n\in\N}$ converges to some $z\in Z$. Coordinatewise inspection shows that $z=y$, so $y\in Z$, and the inequality~\eqref{Def:FIM:eq} follows from~\eqref{L:c00:eq} and continuity of the norms.
\end{proof}  

The final ingredient we need before we can present the proof of \Cref{Thm:main}\ref{thm:main:item2} is a relation~$\preceq$ defined on the family 
\[ \operatorname{BSp} = \{B_p : 1<p<\infty\}\cup\{\ell_p : 1\le p<\infty\}\cup\{ S_p : 1\le p<\infty\}\cup\{c_0\} \] 
of Banach spaces. Its definition is as follows.
\begin{definition}\label{Defn:preceq}
  Let $Y,Z\in\operatorname{BSp}$. Then $Y\preceq Z$ if and only if one of the following five mutually exclusive conditions is satisfied: 
  \begin{itemize}
    \item $Y=\ell_1$ and $Z\in\operatorname{BSp}$ is arbitrary;
    \item $Y=B_p$ for some $1<p<\infty$ and\\ $Z\in\{B_q : p\le q<\infty\}\cup\{\ell_q : 1< q<\infty\}\cup\{ S_q : 1\le q<\infty\}\cup\{c_0\}$;
    \item $Y=\ell_p$ for some $1<p<\infty$ and $Z\in\{\ell_q : p\le q<\infty\}\cup\{ S_q : p\le q<\infty\}\cup\{c_0\}$;
    \item $Y=S_p$ for some $1\le p<\infty$ and $Z\in\{\ell_q : p<q<\infty\}\cup\{ S_q : p\le q<\infty\}\cup\{c_0\}$;
    \item $Y=c_0$ and $Z=c_0$.
  \end{itemize}
  In line with standard practice, we write $Y\prec Z$ when $Y\preceq Z$ and $Y\ne Z$.
\end{definition}

It is easy to see that $\preceq$ is a linear order, whose definition we can summarize as follows:
\begin{equation}\label{Eq:prec}
    \ell_1\prec B_p\prec B_q\prec S_1\prec\ell_p\prec S_p\prec\ell_q\prec S_q\prec c_0\qquad (1<p<q<\infty).
\end{equation}
The next lemma explains its relevance for our purposes.  
\begin{lemma}\label{L:order}
  Let $Y,Z\in\operatorname{BSp}$ with $Y\prec Z$. Then the pair $(Y,Z)$ admits a formal inclusion map \mbox{$R_{Y,Z}\colon Y\to Z$} which is strictly singular.
\end{lemma}

\begin{proof} Admitting a formal inclusion map is clearly a transitive relation in the sense that if the pairs $(X,Y)$ and $(Y,Z)$ both admit formal inclusion maps, then so does the pair $(X,Z)$, and in this case $R_{X,Z} = R_{Y,Z}R_{X,Y}$, which implies that $R_{X,Z}$ is strictly singular whenever at least one of the formal inclusion maps $R_{X,Y}$ and $R_{Y,Z}$ is.  

Hence, in view of~\eqref{Eq:prec}, it suffices to show that each of the following pairs admits a formal inclusion map which is strictly singular: 
\begin{enumerate}[label={\normalfont{(\roman*)}}] 
\item\label{L:order:1} $(\ell_1,B_p)$ for $1<p<\infty$;
\item\label{L:order:2} $(B_p,B_q)$ for $1<p<q<\infty$;
\item\label{L:order:3} $(B_p,S_1)$ for $1<p<\infty$;
\item\label{L:order:4} $(S_p,\ell_q)$ for $1\le p<q<\infty$;
\item\label{L:order:5} $(\ell_p,S_p)$ for $1<p<\infty$;
\item\label{L:order:6} $(S_p,c_0)$ for $1<p<\infty$.
\end{enumerate}

We begin by explaining why these pairs admit formal inclusion maps. By~\Cref{L:c00}, we must show 
that there is a constant $C_{Y,Z}>0$ such that $\lVert x\rVert_Z\le C_{Y,Z}\lVert x\rVert_Y$ for every $x\in c_{00}$.
This is only non-trivial in case~\ref{L:order:4}, so we leave it last. 
In the five other cases, we can easily verify that $C_{Y,Z}=1$ works:
\begin{enumerate}[label={\normalfont{(\roman*)}}] 
\item This is simply subadditivity of the $B_p$-norm together with the fact that the unit vector basis for~$B_p$ is normalized.
\item This is a consequence of the inequality $\beta_q(x,\mathcal{C})\le\beta_p(x,\mathcal{C})$ for every Schreier chain~$\mathcal{C}$, which in turn follows from the well-known fact that the pair $(\ell_p,\ell_q)$ admits a formal inclusion map with constant $C_{\ell_p,\ell_q}=1$.
\item This is immediate from the fact that $\mu_1(x,F) = \beta_p(x,\{F\})$ for every $F\in\mathcal{S}_1$.
\stepcounter{enumi}
\item This is clear because $\mu_p(x,F)\le\lVert x\rVert_{\ell_p}$ for every $F\in\mathcal{S}_1$.
\item This follows from the fact that the coordinate functionals $(e_n^*)_{n\in\N}$ on~$S_p$ have norm~$1$. 
\end{enumerate}

We address case~\ref{L:order:4} by modifying an argument originally due to Graham Jameson for $p=1$, presented in the first part of the proof of \cite[Theorem~A.1]{LS1}, and include the details here for the convenience of the reader. 
Since the $S_p$- and $\ell_q$-norms depend only on the moduli of the coordinates of~$x\in c_{00}$, we may suppose that $x(n)\ge 0$ for every $n\in\N$. We can find a permutation $\sigma\colon\N\to\N$ such that $x\circ\sigma$ is decreasing because~$x$ has finite support, and \cite[Lemma~A.2]{LS1} implies that $\lVert x\circ\sigma\rVert_{S_p}\le \lVert x\rVert_{S_p}$, while $\lVert x\circ\sigma\rVert_{\ell_q} = \lVert x\rVert_{\ell_q}$. Hence, by replacing~$x$ with the decreasing sequence~$x\circ\sigma$, we may suppose that~$x$ is decreasing. Furthermore, by homogeneity, we may suppose that $\lVert x\rVert_{S_p}=1$.  

Take $n\in\N_0$, and set $F_n = [2^n,2^{n+1})\cap\N\in\mathcal{S}_1$. Then $x(2^n)\ge x(j)\ge x(2^{n+1})$ for $j\in F_n$ because $x$ is decreasing. The lower of these bounds implies that
\[ 1 = \lVert x\rVert_{S_p}^p\ge \mu_p(x,F_n)^p=\sum_{j\in F_n} x(j)^p\ge \lvert F_n\rvert\cdot x(2^{n+1})^p = 2^nx(2^{n+1})^p, \]
so $x(2^{n+1})\le 2^{-n/p}$. Consequently we have
\begin{align*} \mu_q(x,F_{n+1})^q &= \sum_{j\in F_{n+1}} x(j)^q = \sum_{j\in F_{n+1}} x(j)^{q-p}x(j)^p\\ &\le x(2^{n+1})^{q-p}\sum_{j\in F_{n+1}} x(j)^p\le  2^{\frac{-n(q-p)}{p}}\mu_p(x,F_{n+1})^p\le \bigl(2^{\frac{p-q}{p}}\bigr)^n, \end{align*}
and therefore
\begin{align*} \lVert x\rVert_{\ell_q}^q &= x(1)^q+\sum_{n=0}^\infty \mu_q(x,F_{n+1})^q\le 1+ \sum_{n=0}^\infty \bigl(2^{\frac{p-q}{p}}\bigr)^n = 1+ \frac{1}{1-2^{\frac{p-q}{p}}} = \frac{2^{\frac{q}{p}}-1}{2^{\frac{q-p}{p}}-1}. \end{align*}    
This proves that the inequality~\eqref{L:c00:eq} is satisfied for some positive constant $C_{S_p,\ell_q}\le \bigl(\frac{2^{q/p}-1}{2^{(q-p)/p}-1}\bigr)^{1/q}$, thus completing the proof in case~\ref{L:order:4}. 

To justify that the formal inclusion maps are strictly singular in each of the six cases above, we require two standard notions concerning a pair of infinite-dimensional Banach spaces $Y$ and~$Z\colon$ 
\begin{itemize}
    \item $Y$ and $Z$ are \emph{totally incomparable} if no infinite-dimensional Banach space embeds isomorphically in both~$Y$ and~$Z$.
    \item $Y$ is \emph{saturated} with copies of~$Z$ if every closed, infinite-dimensional subspace of~$Y$ contains a subspace which is isomorphic to~$Z$. 
\end{itemize}
Trivially, \emph{every} operator between a pair of totally incomparable Banach spaces is strictly singular, and the following three facts imply that each of the first five pairs in the list above are  totally incomparable: 
\begin{itemize}
\item Any pair of distinct spaces from the family $\{\ell_p : 1\le p<\infty\}\cup\{c_0\}$ are totally incomparable, and every space belonging to this family is saturated with copies of itself (see \cite[Proposition~2.a.2, the remark following it, and page~75]{LT1}).
\item $B_p$ is saturated with copies of~$\ell_p$ for $1<p<\infty$ (see \cite[Theorem~II.3.3]{S}, \cite[Theorem~0.15(e)]{CS} or \cite[Theorem~2.4]{LS1}).
\item $S_p$ is saturated with copies of~$c_0$ for $1\le p<\infty$  (see \cite[Corollary~5.4]{BL} or \cite[Theorem~2.4]{LS1}). 
\end{itemize}
Hence, the formal inclusion map $R_{Y,Z}\colon Y\to Z$ is strictly singular in cases~\ref{L:order:1}--\ref{L:order:5}. 

This argument does not work for the formal inclusion map $R_{S_p,c_0}\colon S_p\to c_0$, but \cite[Proposition~6.6]{LS1} shows that it is strictly singular; alternatively, we can easily deduce this result from case~\ref{L:order:4} because $R_{S_p,c_0} = R_{\ell_q,c_0}R_{S_p,\ell_q}$ for $1\le p<q<\infty$. 
\end{proof}

\begin{proof}[Proof of Theorem~{\normalfont{\ref{Thm:main}}}\ref{thm:main:item2}] The family
\begin{equation}\label{Pf:mainthm2:eq1} 
\Sigma = \begin{cases} \{ B_p : p\in L\}\cup\{\ell_q : q\in L\cup M\}\ &\text{if}\ N=\emptyset,\\
\{ B_p : p\in L\}\cup\{\ell_q : q\in L\cup M\}\cup\{S_r : r\in N\}\cup\{c_0\}\ &\text{otherwise} \end{cases} 
\end{equation}
consists of complemented subspaces of the Banach space~$X$ defined by~\eqref{Thm:main:eqX} because~$B_p$ contains a complemented copy of~$\ell_p$ for $1<p<\infty$ and~$S_r$ contains a complemented copy of~$c_0$ for $1\le r<\infty$.

Comparing the definitions~\eqref{Pf:mainthm2:eq1} and~\eqref{Thm:main:eqk}, we see that~$\Sigma$ has cardinality~$k+1$. Since $\Sigma\subset\operatorname{BSp}$, we can use the linear order~$\preceq$ from~\Cref{Defn:preceq} to enumerate its members in increasing order: 
\[ X_{1}\prec X_{2}\prec\cdots\prec X_{k+1}. \]
\Cref{L:order} implies that we have a strictly singular formal inclusion map \mbox{$R_{X_{j},X_{j+1}}\colon X_{j}\to X_{j+1}$} for each $1\le j\le k$. Their composition
\begin{equation*} \spreaddiagramcolumns{5ex}%
     \xymatrix{X_1\ar^-{\displaystyle{R_{X_{1},X_{2}}}}[r] & X_2\ar^-{\displaystyle{R_{X_{2},X_{3}}}}[r] & X_3\ar[r] & \cdots\ar[r] & X_k\ar^-{\displaystyle{R_{X_{k},X_{k+1}}}}[r] & X_{k+1}} \end{equation*}
is simply the formal inclusion map $R_{X_1,X_{k+1}}\colon X_1\to X_{k+1}$, which is not compact because it maps the unit vector basis for~$X_1$ onto the unit vector basis for~$X_{k+1}$. Now the conclusion follows from \Cref{L:complemented}.
\end{proof}

\begin{remark}\label{R:includingc0} The aim of this remark is to justify the two statements made in  \Cref{R:LMdisjoint}\ref{R:LMdisjoint:4}. We recall that $Y=X\oplus c_0$, where~$X$ is given by~\eqref{Thm:main:eqX}. 
\begin{itemize}
    \item Suppose that $N\ne\emptyset$, and take $r\in N$. Arguing as in \Cref{R:LMdisjoint}\ref{R:LMdisjoint2}, we write $S_r\cong W\oplus c_0$ for some Banach space~$W$, which in combination with the fact that $c_0\oplus c_0\cong c_0$ implies that $S_r\oplus c_0\cong W\oplus c_0\oplus c_0\cong W\oplus c_0\cong S_r$, so $Y\cong X$ in this case. 
    \item Now suppose that $N=\emptyset$, so that $Y = \bigl(\bigoplus_{p\in L} B_p\bigr)\oplus\bigl(\bigoplus_{q\in M}\ell_q\bigr)\oplus c_0$ for some finite sets $L\subset(1,\infty)$ and $M\subset[1,\infty)$, and set $k=\lvert L\rvert + \lvert L\cup M\rvert$. The proof of \Cref{Thm:main}\ref{thm:main:item1} that we gave above in the case $N=\emptyset$ carries over almost verbatim because \Cref{L:compactnessofSScompositions}\ref{L:compactnessofSScompositions1} applies to~$c_0$, too; that is, every strictly singular operator on~$c_0$ is compact. 

    It is also easy to modify the above proof of \Cref{Thm:main}\ref{thm:main:item2} to the present context: Simply define $\Sigma = \{ B_p : p\in L\}\cup\{\ell_q : q\in L\cup M\}\cup\{c_0\}$ and argue as before. 
\end{itemize}    
\end{remark}

According to Pitt's Theorem \cite[Proposition~2.c.3]{LT1}, every operator from~$\ell_p$ to~$\ell_q$ is compact for $1\le q<p<\infty$. Seifert claimed in his dissertation \cite[Corollary~II.3.4]{S} that the analogous result is true for the Baernstein spaces, that is, every operator from~$B_p$ to~$B_q$ is compact for $1<q<p<\infty$. However, using \Cref{L:order}, we can easily show that Seifert's claim is false. The origin of this error appears to be \cite[Lemma~II.3.2]{S}; although Seifert does not explicitly cite this lemma in his proof of \cite[Corollary~II.3.4]{S}, he uses it implicitly. We refer to \cite[the paragraph below Theorem~2.4]{LS1} for a corrected version of \cite[Lemma~II.3.2]{S}. 

Unfortunately, Seifert's incorrect claim has gained much wider publicity than one would usually expect for a result contained in an unpublished PhD thesis because it was reproduced (without proof) in the lecture notes~\cite[Theorem~0.15(f)]{CS}. Hence, we shall state our correction formally. 

\begin{corollary}\label{C:Seifertmistake}
    There are strictly singular, non-compact operators from~$B_p$ to~$B_q$ for every pair $p,q\in(1,\infty)$. 
\end{corollary}

\begin{proof} \Cref{Defn:preceq} shows that $B_p\prec\ell_q$, so by \Cref{L:order}, we have a formal inclusion map $R_{B_p,\ell_q}\colon B_p\to\ell_q$ which is strictly singular. It is not compact because it maps the unit vector basis for~$B_p$ onto the unit vector basis for~$\ell_q$. Since~$B_q$ contains a complemented subspace that is isomorphic to~$\ell_q$, we can find operators $U\colon B_q\to\ell_q$ and $V\colon \ell_q\to B_q$ such that $I_{\ell_q} = UV$. It follows that $VR_{B_p,\ell_q}\colon B_p\to B_q$ is a  strictly singular, non-compact operator.
\end{proof}

The realization that there are non-compact operators from~$B_p$ to~$B_q$ for $1<q<p<\infty$ raises the question whether some, possibly very rapidly increasing, subsequence of the unit vector basis for~$B_p$ dominates a subsequence of the unit vector basis for~$B_q$. We conclude by showing that this is impossible.

\begin{proposition} Let $p,q\in(1,\infty)$. Then the unit vector basis for~$B_p$ admits a subsequence which dominates a subsequence of the unit vector basis for~$B_q$ if and only if $p\le q$.     
\end{proposition}

\begin{proof} For clarity, we denote the unit vector bases for~$B_p$ and~$B_q$ by $(e_n^p)_{n\in\N}$ and $(e_n^q)_{n\in\N}$, respectively. 

    The implication $\Leftarrow$ is clear because, for $1<p\le q<\infty$, the formal inclusion map $B_p\to B_q$ is bounded, which means that $(e_n^p)_{n\in\N}$ dominates $(e_n^q)_{n\in\N}$. 

    Conversely, suppose that $(e_{m_j}^p)_{j\in\N}$ dominates $(e_{n_j}^q)_{j\in\N}$ for some integers $1\le m_1<m_2<\cdots$ and $1\le n_1<n_2<\cdots$. Then  
    \[ x_k = \frac1{2^{k-1}}\sum_{j=2^{k-1}}^{2^k-1} e_{m_j}^p\in B_p\qquad\text{and}\qquad y_k = \frac1{2^{k-1}}\sum_{j=2^{k-1}}^{2^k-1} e_{n_j}^q\in B_q\qquad (k\in\N) \]
    are unit vectors because their supports $\{ m_j : 2^{k-1}\le j<2^k\}$ and $\{ n_j : 2^{k-1}\le j<2^k\}$ are Schreier sets. By hypothesis, the block basic sequence $(x_k)_{k\in\N}$ dominates $(y_k)_{k\in\N}$. Since $\lVert x_k\rVert_\infty = 1/2^{k-1}\to0$ as $k\to\infty$, \cite[Proposition~2.14]{LS1} implies that $(x_k)_{k\in\N}$ admits a subsequence $(x_{k_j})_{j\in\N}$ which is dominated by the unit vector basis for~$\ell_p$. Furthermore, being a normalized block basic sequence of the unit vector basis for~$B_q$, $(y_{k_j})_{j\in\N}$ dominates the unit vector basis for~$\ell_q$ by \cite[Lemma~2.10]{LS1}. In conclusion, it follows that the unit vector basis for~$\ell_p$ dominates the unit vector basis for~$\ell_q$, which is possible only if $p\le q$. 
\end{proof}

\begin{acknowledgements} The bulk of the research on which this paper is based was carried out during a research retreat at the Isaac Newton Institute (INI) in Cambridge, UK, in May 2025. We are grateful to the INI for their kind hospitality and financial support that made this visit possible. Wirzenius is supported by project L100192451 of the Czech Academy of Sciences. His visit to the UK received additional financial support from a grant of the Ruth and Nils-Erik Sten\-b\"{a}ck Foundation. He acknowledges these sources of funding with thanks.

For the purpose of open access, the authors have applied a Creative Commons Attribution (CC-BY) licence to any Author Accepted Manuscript version arising.
\end{acknowledgements}

\end{document}